\documentclass[12pt]{amsart}
\usepackage[alphabetic]{amsrefs}
\usepackage[bookmarksopen,bookmarksdepth=2,breaklinks=true]{hyperref}
\usepackage[hyphenbreaks]{breakurl}
\usepackage{verbatim}
\usepackage{comment}

\def\spp{\vspace{5pt} \noindent}

\usepackage[margin=1.5cm]{geometry}                
\usepackage{graphicx}
\usepackage{xcolor}

\usepackage{amsmath, mathtools, amsfonts, amssymb}

\numberwithin{equation}{section}

\newtheorem{theorem}{Theorem}[section]
\newtheorem{example}[theorem]{Example}

\newtheorem{definition}[theorem]{Definition}
\newtheorem{remark}[theorem]{Remark}

\renewcommand{\le}{\leqslant}
\renewcommand{\leq}{\leqslant}
\renewcommand{\ge}{\geqslant} 
\renewcommand{\geq}{\geqslant}

\allowdisplaybreaks

\author{Serena Dipierro}
\author{Lyle Noakes}
\author{Enrico Valdinoci}

\address{Department of Mathematics and Statistics,
University of Western Australia, 35 Stirling Highway,
WA6009 Crawley, Australia}
\email{serena.dipierro@uwa.edu.au, lyle.noakes@uwa.edu.au, enrico.valdinoci@uwa.edu.au}

\title{Napoleonic triangles on the sphere}

\begin{document}
\maketitle

\begin{abstract} As is well-known, numerical experiments show that Napoleon's Theorem for planar triangles does not extend to a similar statement for triangles on the unit sphere $S^2$. Spherical triangles for which an extension of Napoleon's Theorem holds are called {\em Napoleonic}, and until now the only known examples have been equilateral. 
In this paper we determine all Napoleonic spherical triangles, including a class corresponding to points on a 2-dimensional ellipsoid, whose Napoleonisations are all congruent. Other new classes of examples are also found, according to different versions of Napoleon's Theorem for the sphere. The classification follows from successive simplifications of a complicated original algebraic condition, exploiting geometric symmetries and algebraic factorisations. \end{abstract}

\section{Introduction}

\subsection{Napoleonic questions}\label{Nq}
A beautiful, theorem, often attributed to
Napoleon Bonaparte, is that
if equilateral triangles are constructed either all outwards or all inwards on the sides of a triangle in the Euclidean plane $E^2$, then the lines connecting the centres of those equilateral triangles themselves form an equilateral triangle. \spp

\begin{theorem}[Napoleon's Theorem]\label{NAPTHO}
Given distinct ~$P_0$, $P_1$, $P_2\in E^2$, let~$Q_0$, $Q_1$, $Q_2\in E^2$ be such that, for every~$i\in \mathbb{Z}/3\mathbb{Z}$, 
the triangle $Q_iP_{i+1}P_{i+2}$ is
equilateral with interior disjoint from the interior of the triangle~$P_0P_1P_2$. Then $R_0R_1R_2$ is equilateral, where $R_i$ is the centroid of $Q_iP_{i+1}P_{i+2}$. 

\spp
If instead ``disjoint''  is replaced by ``not disjoint'' the triangle $R_0R_1R_2$  is also equilateral.~$\square$
\end{theorem}

\spp
We refer to ~\cite{MR635364, MR2928662} for some of the history of Napoleon's Theorem, including discussion of  its attribution. For    
 extensions of Napoleon's Theorem see ~\cite{MR1847491, MR2410581, MR2558297}
and the references therein
(also {\tt https://www.cut-the-knot.org/proofs/napoleon.shtml}
for a collection of related work). Even in flat space, these kinds of elementary questions have attracted continued attention. For example, the PDN Theorem ~\cite{PETR-zbMATH02641104, DOUG-MR2178, NEUM-MR6839}, extending Napoleon's Theorem from planar triangles to polygons, was found independently by three first-rate mathematicians, including Fields Medallist Jesse Douglas. 

\spp
Our extension of Theorem \ref{NAPTHO} is a classification of Napoleonic triangles in the unit $2$-sphere $S^2$ in Euclidean $3$-pace $E^3$. Here vertices are points in $S^2$, and the distance  $d(P,Q)\in [0,\pi ]$ between $P,Q\in S^2$ is $\arccos \langle P,Q\rangle$, where $\langle ~,~\rangle$ is the Euclidean inner product. Edges of spherical triangles are defined to be minimal great circle arcs $PQ$. Then for vertices $P_0,P_1,P_2\in S^2$ where $P_i\not= \pm P_{i+1}$, the {\em spherical triangle} $P_0P_1P_2$ is a union of edges $P_iP_{i+1}$. It is convenient to consider the indices $i=0,1,2$ as elements of $\mathbb{Z}/3\mathbb{Z}$ namely the integers reduced mod $3$. 

\spp
A well-known quick proof of the planar Napoleon's Theorem uses symmetries and dilations, but this breaks down in spherical geometry, because   
 most spherical triangles are not {\em Napoleonic} in the following sense.
\spp

\begin{definition}\label{NAPTHD}
A {\em Napoleonisation} of a spherical triangle $P_0P_1P_2$ is a spherical triangle of the form $R_0R_1R_2$, 
where, for all $i\in \mathbb{Z}/3\mathbb{Z}$, $R_i$ is the centroid of an equilateral spherical triangle $Q_iP_{i+1}P_{i+2}$.  We say $P_0P_1P_2$ is {\em Napoleonic} when it has an equilateral Napoleonisation.~$\square$
\end{definition}

\clearpage

\spp
For there to exist a vertex $Q_i\in S^2$ with $Q_iP_{i+1}P_{i+2}$ equilateral, we must have $d(P_{i+1},P_{i+2})\leq 2\pi /3$ for all $i\in \mathbb{Z}/3\mathbb{Z}$. So assume from now on that $d(P_{i+1},P_{i+2})\leq 2\pi /3$. Notice also that
\begin{itemize}
\item when $d(P_{i+1},P_{i+2})=2\pi /3$ the vertex $Q_i$ is unique, and $Q_iP_{i+1}P_{i+2}$ are cogeodesic, namely they lie on a great circle. 
\item when $0<d(P_{i+1},P_{i+2})<2\pi /3$ there are precisely two choices of $Q_i$.
\end{itemize}

\spp
\begin{example}\rm Suppose that $P_0,P_1,P_2$ are cogeodesic: then they are all orthogonal to 
some $N\in S^2$, and 
$$\sum_{i=0}^2d(P_i,P_{i+1})=2\pi \Longrightarrow d(P_i,P_{i+1})=2\pi /3\hbox{ for all }i\in \mathbb{Z}/3\mathbb{Z}.$$
So $P_0P_1P_2$ is already equilateral and, for all $i\in \mathbb{Z}/3\mathbb{Z}$, we have $Q_i=P_i$ and $R_i=\pm N$. 
So the $R_i$ are not distinct, $R_0R_1R_2$ is not considered as a triangle, and $P_0P_1P_2$ is not Napoleonic. 
%
%
On the other hand, if $P_0P_1P_2$ is equilateral, with $d(P_i,P_{i+1})<2\pi /3$ then, by symmetry, $P_0P_1P_2$ is Napoleonic.  $\square$
\end{example}

\spp
Suppose from now on that, in addition to $d(P_i,P_{i+1})\leq 2\pi /3$, 
the $P_i$ are distinct and not cogeodesic. The {\em interior} of the spherical triangle $P_0P_1P_2$ is the path-component of the spherical centroid (or barycentre)
$$\frac{P_0+P_1+P_2}{\Vert P_0+P_1+P_2\Vert}$$
in $S^2\setminus P_0P_1P_2$. 

\spp
The Napoleonisation $R_0R_1R_2$ is said to be {\em outward} (respectively {\em inward}) when all $Q_i$ are directed away from (respectively towards) the interior of the spherical triangle $P_0P_1P_2$.
More detail is given in \S\ref{SEZ12} below. 
We call $P_0P_1P_2$ {\em outward-Napoleonic} or {\em inward-Napoleonic}  when its outward Napoleonisation or inward Napoleonisation is equilateral, respectively. To begin we consider outward Napoleonisations. 
 
\begin{figure}[h]
\centering
\includegraphics[width=6cm]{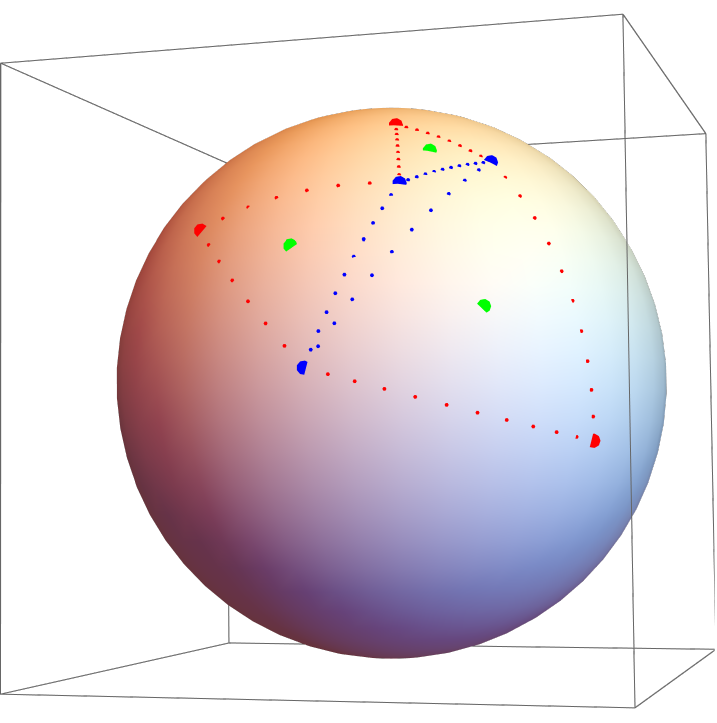}
\caption{\em Vertices (green) of the outward Napoleonisation of the spherical triangle~$P_0P_1P_2$, with~$P_0:=(1,0,0)$,
$P_1:=\left( \frac34,\frac14,\frac{\sqrt6}{4}\right)$ and~$P_2:=\left( \frac12,\frac12,\frac{\sqrt2}{2}\right)$ (blue).}
\label{FESEedIG1EPS32}\end{figure}

\spp
Figure~\ref{FESEedIG1EPS32} shows the vertices (green) $R_0,R_1,R_2$ of the outward Napoleonisation for a particular choice of $P_0,P_1,P_2$ in the unit sphere $S^2$. The spherical distances $d(R_0,R_1),d(R_1,R_2),d(R_2,R_0)$ are 
$$0.474066,0.490357,0.448728.$$ 
So $R_0R_1R_2$ is not equilateral, namely $P_0P_1P_2$ is not outward-Napoleonic. 

\spp
Next we give some explicit examples of Napoleonic constructions on the sphere.

\subsection{Napoleonic examples and counterexamples}\label{SEZ11}
As noted in ~\cite{TREUDEN} and again in \S\ref{Nq}, examples of non-Napoleonic spherical triangles are easily found using Mathematica. Examples can also be found using mathematical reasoning as follows. 

\begin{example}\label{Nn}{\rm
Let~$\delta>0$, to be taken as small as we wish.
We denote by~$o(1)$ quantities which are infinitesimal as~$\delta\searrow0$.

\spp
Let~$P_0:=(1,0,0)$ and~$P_1:=\left( -\frac12,\frac{\sqrt3}2,0\right)$.
Let also~$P_2\in S^2$ be such that
$$ \Vert P_0-P_2\Vert<\delta\qquad {\mbox{ and }}\qquad \langle P_2,\,(0,1,0)\rangle>0.$$

\spp
Notice that there is one, and only one, possible choice of~$Q_2\in S^2$ such that the spherical triangle~$P_0P_1Q_2$ is equilateral, namely~$Q_2:=\left( -\frac12,-\frac{\sqrt3}2,0\right)$ (which is actually the largest possible equilateral spherical triangle in~$S^2$). There are however two spherical triangles~$P_0P_1Q_2$, one containing~$(0,0,1)$ and one containing~$(0,0,-1)$. Let us pick the first one (the other choice leading to similar calculations). In this case, the centroid~$R_2$ (taken as the intersection of the medians) of the spherical triangle~$P_0P_1Q_2$ is~$(0,0,1)$.

\spp
We also observe that the spherical distance~$\vartheta_0$ between~$P_1$ and~$P_2$ is such that
$$\cos\vartheta_0=\langle P_1,P_2\rangle
=\left\langle \left( -\frac12,\frac{\sqrt3}2,0\right),\,P_2\right\rangle
>\left\langle \left( -\frac{1}2,0,0\right),\,P_2\right\rangle
\ge-\frac{1}2,$$ therefore there are two possible choices for~$Q_0\in S^2$ such that the spherical triangle~$P_1P_2Q_0$ is equilateral. Anyway, since the choice of~$Q_2$ was unique, $Q_0$ needs to be at an~$o(1)$-distance from~$Q_2$.
Consequently, the centroid~$R_0$ of~$P_1P_2Q_0$ is also~$o(1)$-close to
either~$(0,0,1)$ or~$(0,0,-1)$.

\spp
Also, since~$P_0$ and~$P_2$ lie at an~$o(1)$-distance from each other, either equilateral spherical triangle~$P_0P_2Q_1$ (there are two
possible choices of~$Q_1\in S^2$) is contained in an~$o(1)$-neighborhood of~$P_0$, therefore its centroid~$R_1$ is also contained in a~$o(1)$-neighborhood of~$P_0$.

\spp
{F}rom these considerations, 
when~$R_0$ is~$o(1)$-close to~$(0,0,1)=R_2$,
for small~$\delta$ we arrive at
$$ \Vert R_1-R_2\Vert
\ge\Vert P_0-R_2\Vert-o(1)
=\Vert (1,0,0)-(0,0,1)\Vert-o(1)=\sqrt2-o(1)>o(1)=\Vert R_0-R_2\Vert,$$
while, when~$R_0$ is~$o(1)$-close to~$(0,0,-1)$, we see that
\begin{eqnarray*}&& \Vert R_1-R_2\Vert
\le\Vert P_0-R_2\Vert+o(1)
=\Vert (1,0,0)-(0,0,1)\Vert+o(1)=\sqrt2+o(1)<2 \\&&\qquad=
\Vert (1,0,0)-(-1,0,0)\Vert
=\Vert R_0-R_2\Vert,\end{eqnarray*}
showing that the spherical triangle~$R_0R_1R_2$ is not equilateral.~$\square$}\end{example}

\begin{remark}\label{Rou}{\rm
Example~\ref{Nn} can be also understood by perturbing the degenerate triangle $P_0P_1P_0$ where 
 ~$P_0:=(1,0,0)=:P_2$ and~$P_1:=\left( -\frac12,\frac{\sqrt3}2,0\right)$. In this case~$Q_1=(1,0,0)$ and~$Q_0=Q_2=\left( -\frac12,\frac{\sqrt3}2,0\right)$. Furthermore, $R_1=(1,0,0)$ and~$R_0$, $R_2\in\{ (0,0,1),(0,0,-1)\}$.
We thus obtain that~$\Vert R_1-R_2\Vert= \sqrt2$
while~$\Vert R_0-R_2\Vert\in \{ 0,2\}$, showing that the spherical triangle~$R_0R_1R_2$ is not equilateral. Although $P_0P_1P_0$ does not have $3$ distinct vertices, this defect is removed by suitably small perturbations.~$\square$}
\end{remark}

\spp
By symmetry, both the outward and inward Napoleonisations of an equilateral triangle are also equilateral, with $R_0=R_1=R_2$ for the inward Napoleonisation. So equilateral spherical triangles are both outward-Napoleonic and inward-Napoleonic. 

\spp
Some Napoleonic spherical triangles are not equilateral, indeed not even isosceles. 
\begin{example} \label{NAPEEX} {\rm
Let
$$ P_0 := (1, 0, 0),\qquad P_1:=\left(-\frac7{50},\frac{\sqrt{2451}}{50}, 0\right)\qquad{\mbox{and}}\qquad
P_2: = \left(-\frac{17}{50}, -\frac{423 \sqrt{2451}}{40850},\frac{ 46 \sqrt{4902}}{4085}\right).$$
Then ~$P_0P_1P_2$~ is not isosceles, because
\begin{equation}\label{LJSD:PSL-NO0} \langle P_0,P_1\rangle=-\frac{7}{50},\qquad\langle P_0,P_2\rangle=-\frac{17}{50}\qquad{\mbox{and}}\qquad\langle P_1,P_2\rangle=-\frac{23}{50}.\end{equation}
Define $Q_0,Q_1,Q_2\in S^2$ by 
 \begin{equation}\label{LJSD:PSL-NO1}\begin{split}&
Q_0:=\left( 0, -\frac{23}{\sqrt{2451}}, -31 \sqrt{\frac{2}{2451}}\right)
,\\& Q_1:=\left( 
-\frac{17}{50}, -\frac{1019}{50 \sqrt{2451}}, -\frac{148}{5}\sqrt{\frac{2}{2451}}
\right)\\
{\mbox{and}}\qquad& Q_2:=\left( 
-\frac{7}{50}, -\frac{7}{50} \sqrt{\frac{57}{43}}, -\frac{3}{5} \sqrt{\frac{114}{43}}.
\right).\end{split}\end{equation}
Then
\begin{eqnarray*}
&& \langle P_0,P_1\rangle=\langle P_0,Q_2\rangle=\langle P_1,Q_2\rangle=-\frac7{50},\\
&& \langle P_0,P_2\rangle=\langle P_0,Q_1\rangle=\langle P_2,Q_1\rangle=-\frac{17}{50}\\ {\mbox{and }}
&& \langle P_1,P_2\rangle=\langle P_1,Q_0\rangle=\langle P_2,Q_0\rangle=-\frac{23}{50},
\end{eqnarray*}
and so ~$P_0P_1Q_2$, $P_0P_2Q_1$ and~$P_1P_2Q_0$ are equilateral, with centroids respectively 
%
\begin{eqnarray*}&& R_2:=\frac{P_0+P_1+Q_2}{\|P_0+P_1+Q_2\|}=\left(
\frac{\sqrt{6}}{5}, \frac{3}{5} \sqrt{\frac{38}{43}}, -\sqrt{\frac{19}{43}}
\right),\\&&
R_1:=\frac{P_0+P_2+Q_1}{\|P_0+P_2+Q_1\|}=\left(
\frac{2}{5} \sqrt{\frac{2}{3}}, -\frac{286}{15} \sqrt{\frac{2}{817}}, -\frac{5}{3 \sqrt{817}}
\right)
\\
{\mbox{and}}\qquad&&
R_0:=\frac{P_1+P_2+Q_0}{\| P_1+P_2+Q_0\|}=
\left(-\frac{2}{5} \sqrt{6}, \frac{8}{15}\sqrt{\frac{2}{817}}, -\frac{17}{3\sqrt{817}}
\right) .
\end{eqnarray*}
As a result,
\begin{equation}\label{JSDM:0qowdjfev} \langle R_0,R_1\rangle=\langle R_0,R_2\rangle=\langle R_1,R_2\rangle=-\frac13.\end{equation}
This shows that~$P_0P_1P_2$ is outward-Napoleonic and not isosceles, as illustrated in 
 Figure~\ref{FESEUNMedIG1EPS32}.$\square$

\begin{figure}[h]
\centering
\includegraphics[width=3.4cm]{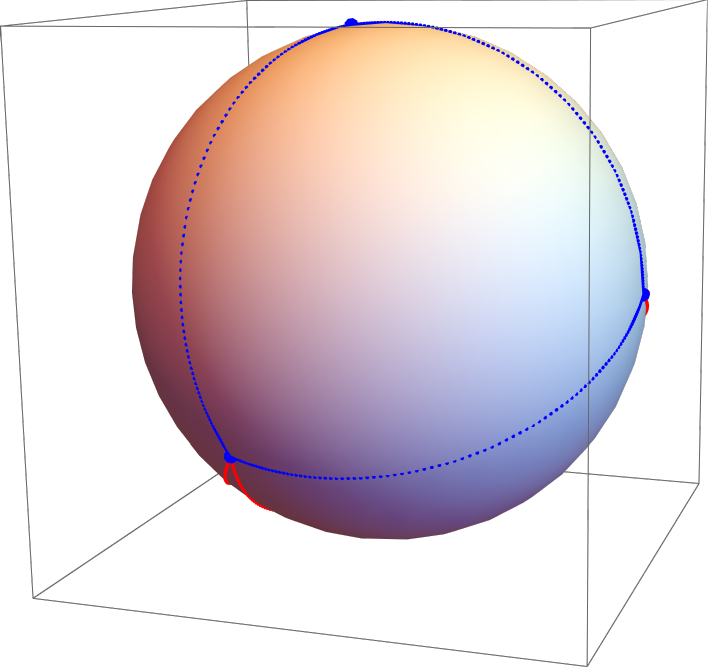}$\qquad$
\includegraphics[width=3.4cm]{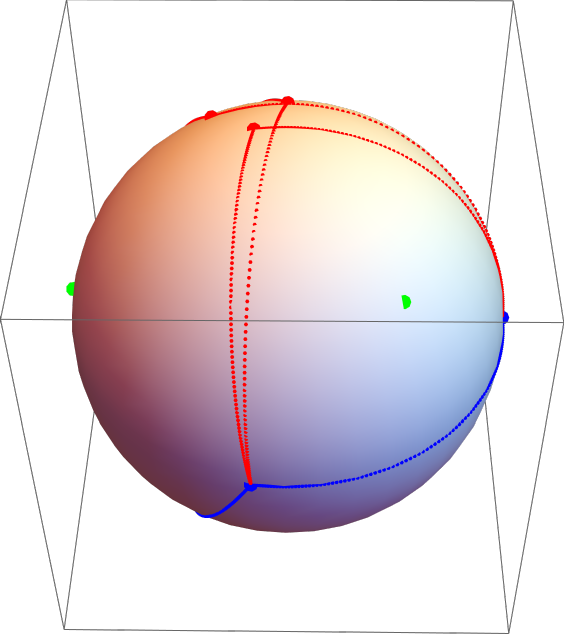}$\qquad$\includegraphics[width=3.4cm]{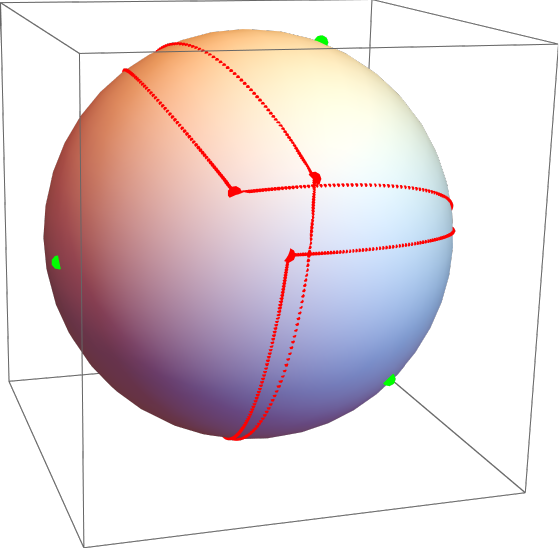}$\qquad$\includegraphics[width=3.4cm]{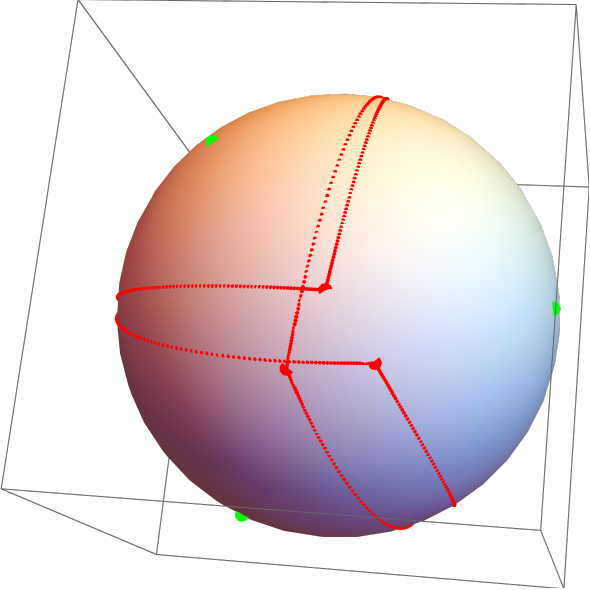}
\caption{\em Different views of the Napoleonic construction related to~\eqref{LJSD:PSL-NO0}
and~\eqref{LJSD:PSL-NO1} (the triangle $P_0P_1P_2$ being in blue and the points~$R_0$, $R_1$, $R_2$ in green).}
\label{FESEUNMedIG1EPS32}\end{figure}

\spp
Unlike the planar case, ~$R_0R_1R_2$ does not always have the same barycentre as $P_0P_1P_2$.
In the present instance, the barycentre of~$P_0P_1P_2$ is
$$\left(\frac{13}{10\sqrt{7}},\frac{197}{10}\sqrt{\frac{3}{5719}},
23\sqrt{\frac{6}{5719}}
\right)\approx
(0.491354,0.451198,0.744978),$$
while the barycentre of~$R_0R_1R_2$ is
$$\left(
-\frac15\sqrt{\frac23},-\frac{107}{15}\sqrt{\frac{2}{817}}, -\frac{79}{3\sqrt{817}}
\right)\approx
(-0.163299, -0.352936, -0.921287).$$
\hfill$\square$}\end{example}

\spp
These examples show that:
\begin{itemize}
\item unlike planar triangles, not all spherical triangles are Napoleonic,
\item as for planar triangles, all equilateral spherical triangles are Napoleonic,
\item there exist non-isosceles outward-Napoleonic spherical triangles,
\item the barycentre of a spherical Napoleonisation is not necessarily the same as the barycentre of the original triangle, even when the latter is Napoleonic.\end{itemize}

\subsection{Classification of Outward and Inward Napoleonic Spherical Triangles}\label{SEZ12} 
For $i\in \mathbb{Z}/3\mathbb{Z}$ let ~$P_i\in S^2$ be distinct, satisfying\footnote{The condition $\langle P_{i+1},P_{i+2}\rangle \geq -1/2$ is needed to ensure that there exist equilateral spherical triangles $Q_iP_{i+1}P_{i+2}$, and $\langle P_{i+1},P_{i+2}\rangle \not= -1/2$ ensures that the interiors of the equilateral triangles are well-defined (otherwise $Q_i,P_{i+1},P_{i+2}$ would be cogeodesic). } 
\begin{equation}\label{COND:0}
\langle P_{i+1},P_{i+2}\rangle > -1/2.\end{equation}
Define 
\begin{equation}\label{DKDEF} d_i:=\sqrt{1+2\langle P_{i+1},P_{i+2}\rangle}\in (0,\sqrt{3}),\end{equation}
determined by the spherical distance between~$P_{i+1}$ and~$P_{i+2}$. Note that 
$\chi :=\langle P_i,P_{i+1}\times P_{i+2}\rangle$ is nonzero and independent of $i\in \mathbb{Z}/3\mathbb{Z}$. 
We relabel the $P_i$ if necessary, to ensure that $\chi >0$. 
Then in \S \ref{SECTION:2} we see that $Q_i$ is directed towards the interior of $P_0P_1P_2$
 precisely when 
$\langle Q_i,P_{i+1}\times P_{i+2}\rangle >0$ (otherwise outwards).

\spp
By symmetry, equilateral spherical triangles are both outward-Napoleonic and inward-Napoleonic. Our main result is that these are the only inward-Napoleonic spherical triangles, and that the non-equilateral outward-Napoleonic spherical triangles are given by a quadratic condition (with all outward-Napoleonisations congruent).

\begin{theorem}\label{THEOREM1}
Suppose that the spherical triangle~$P_0P_1P_2$ is not equilateral. Then $P_0P_1P_2$ is not inward-Napoleonic, and   is outward-Napoleonic if and only if 
\begin{equation}\label{COND-D}
d_0^2+d_1^2+d_2^2+d_0d_1+d_0d_2+d_1d_2=2.
\end{equation}
In such a case, the equilateral Napoleonisation $R_0R_1R_2$ has side $\pi -\arccos (1/3)$.  
\end{theorem}

\begin{figure}[h]
\centering
\includegraphics[width=7cm]{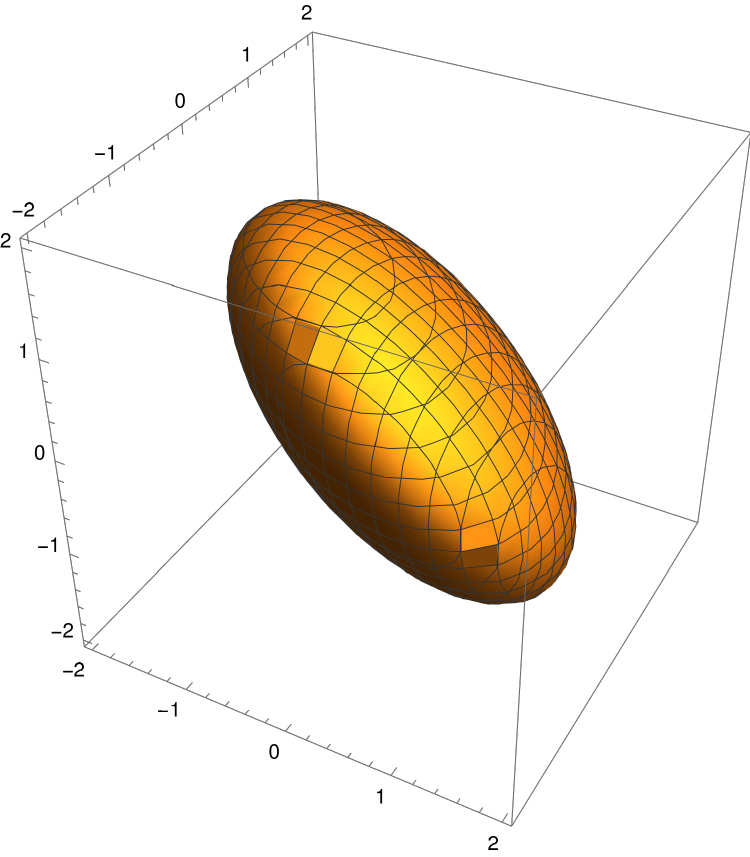}$\qquad$\includegraphics[width=7cm]{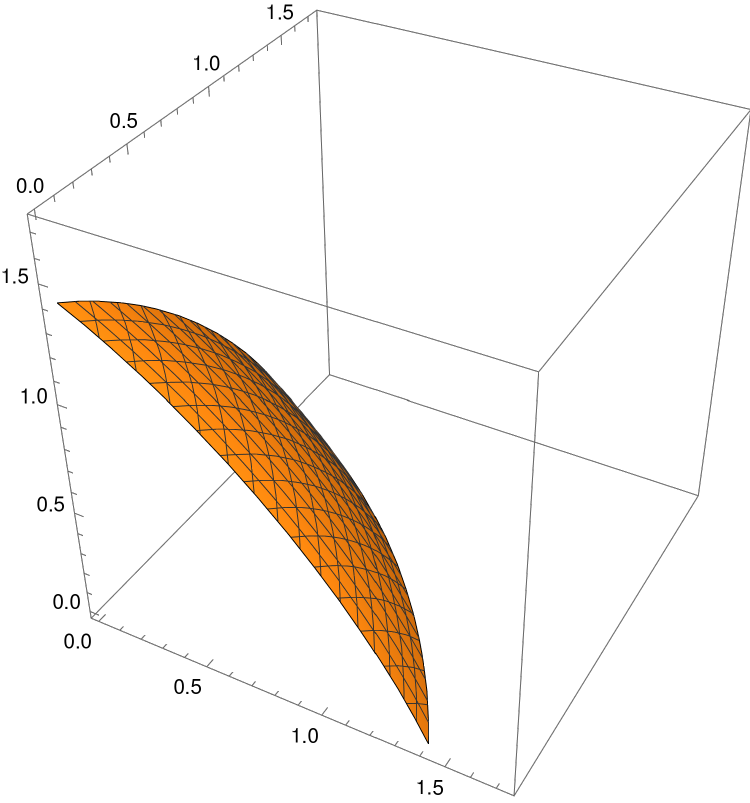}
\caption{\em The ellipsoid of revolution~$\{d_0^2+d_1^2+d_2^2+d_0d_1+d_0d_2+d_1d_2=2\}$ and its intersection
with the region where the $d_i\in (0,\sqrt{3})$.
}
\label{FIG1con}\end{figure}

\spp
So condition~\eqref{COND-D}
is satisfied by $P_0P_1P_2$  \footnote{Indeed,
in Example~\ref{NAPEEX} one has~$d_0=\frac{\sqrt{2}}{5}$,
$d_1=\frac{2\sqrt{2}}{5}$ and~$d_2=\frac{3\sqrt{2}}{5}$, which satisfy condition~\eqref{COND-D}.} in Example~\ref{NAPEEX}, where the dimensions \eqref{JSDM:0qowdjfev} of the Napoleonisation  are also given by Theorem \ref{THEOREM1}. Condition~\eqref{COND-D} is illustrated in Figure~\ref{FIG1con} and, after 
rotating coordinates
$$ \left(\begin{matrix}
X \\ Y \\ Z
\end{matrix}\right)=
\left(\begin{matrix}
1/\sqrt{3} & 1/\sqrt{3} & 1/\sqrt{3} \\
-\sqrt{2/3} & 1/\sqrt{6} & 1/\sqrt{6} \\
0 & -1/\sqrt{2} & 1/\sqrt{2}
\end{matrix}\right)\left(\begin{matrix}
d_0\\ d_1 \\ d_2
\end{matrix}\right)
$$
defines part of\footnote{The conditions $d_i\in (0,\sqrt{3})$ must still be satisfied.} a standard ellipsoid of revolution
$$ 2 X^2 + \frac{Y^2}2 + \frac{Z^2}2=2.$$

\spp
The paper is organised as follows.  
In \S\ref{SECTION:2} we introduce a mathematical framework for studying Napoleonic spherical triangles. 
This allows us to make in \S\ref{SECTION:3}
the arguments required for the proof of Theorem \ref{THEOREM1}, namely the classification of Napoleonic triangles in~$S^2$.
This completes the proofs of our mathematical results. 
In \S\ref{contextsec} we make some more general contextual remarks. Finally in \S\ref{CONC} we list some conclusions.

\section{Spherical Triangles}\label{SECTION:2}
For $P_0\not= \pm P_1\in S^2$ and any $Q\in S^2$ write 
$Q=a_0P_0+a_1P_1+bP_0\times P_1$. 
Then 
\begin{equation}\label{eq1}1=\langle Q,Q\rangle = a_0^2+a_1^2+2a_0a_1\langle P_0,P_1\rangle +b^2(1-\langle P_0,P_1\rangle ^2).\end{equation}
From $a_0+a_1\langle P_0,P_1\rangle =\langle Q,P_0\rangle $ and 
$a_0\langle P_0,P_1\rangle +a_1=\langle Q,P_1\rangle $ 
we obtain 
\begin{equation}\label{eq20}a_0=\frac{\langle Q,P_0\rangle -\langle P_0,P_1\rangle \langle Q,P_1\rangle }{1-\langle P_0,P_1\rangle ^2},~a_1=\frac{\langle Q,P_1\rangle-\langle P_0,P_1\rangle \langle Q,P_0\rangle}{1-\langle P_0,P_1\rangle ^2}.\end{equation} 
Let $P_0QP_1$ be isosceles with $\langle P_0,Q\rangle =\langle P_1,Q\rangle$. Then  (\ref{eq20}) reduces to 
\begin{equation}\label{eq2}
a_0=a_1=a:=\frac{\langle P_0,Q\rangle}{1+\langle P_0,P_1\rangle}.
\end{equation}

\begin{figure}[h]
\centering
\includegraphics[width=5cm]{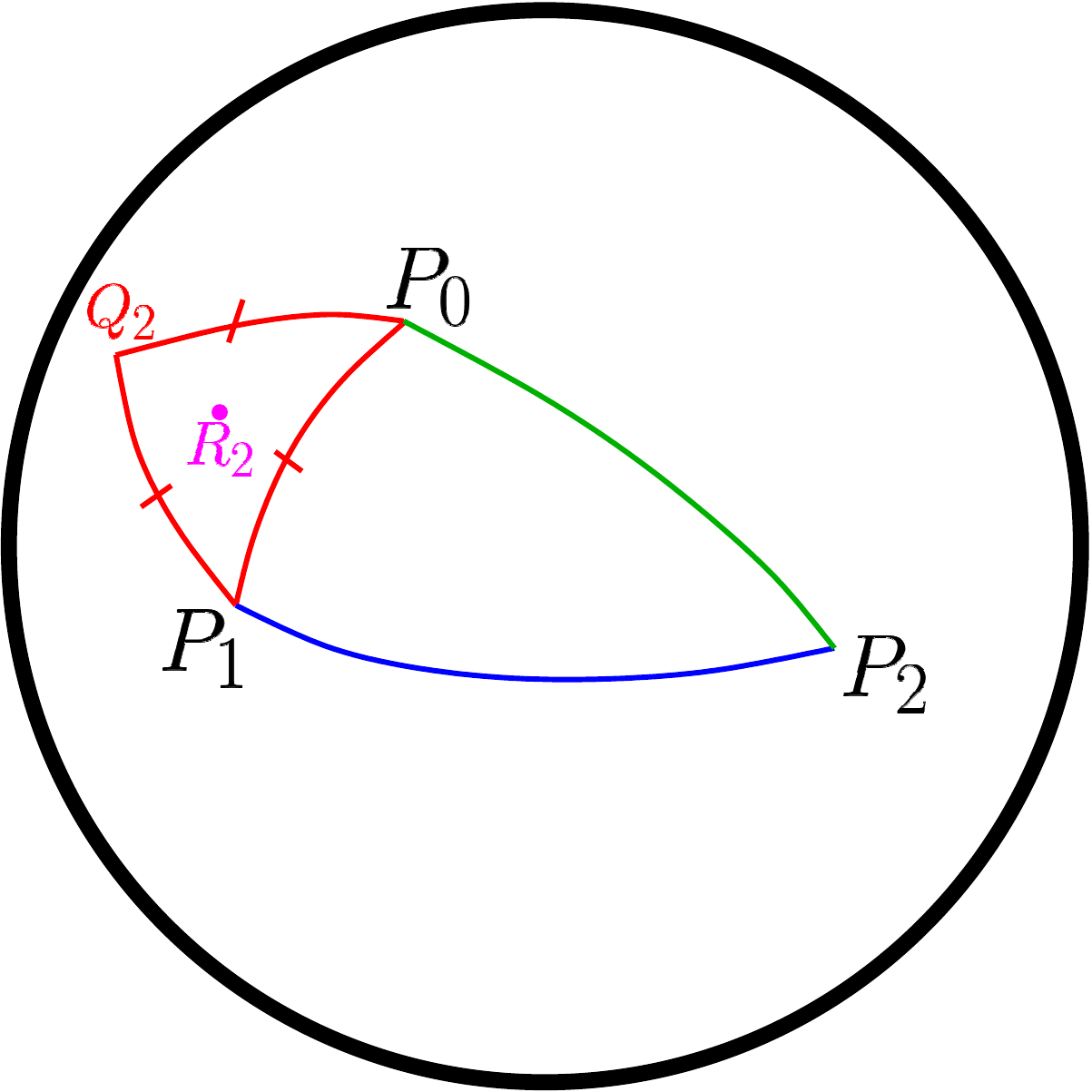}
\caption{\em Constructing an equilateral spherical triangle~$P_0P_1Q_2$, given the spherical triangle~$P_0P_1P_2$.}
\label{FIG1}\end{figure}

\spp
Now take $P_0P_1Q_2$ to be equilateral. 
By (\ref{eq1}), (\ref{eq2}), \begin{equation}\label{eq3}b=\epsilon _2 \frac{\sqrt{1+2\langle P_0,P_1\rangle }}{1+\langle P_0,P_1\rangle }\end{equation}
where $\langle P_0,P_1\rangle \geq -1/2$ in order for the equilateral triangle\footnote{When $\langle P_0,P_1\rangle =-1/2$ the equilateral spherical triangle $P_0P_1Q_2$ is the great circle through $P_0,P_1$.} to exist, and $\epsilon _2=\pm 1$. Then  
\begin{equation}\label{formq2}Q_2=\frac{1}{1+\langle P_0,P_1\rangle}(\langle P_0,P_1\rangle(P_0+P_1)+\epsilon _2\sqrt{1+2\langle P_0,P_1\rangle }P_0\times P_1).\end{equation}
Suppose now that $P_0P_1$ is an edge of a spherical triangle $P_0P_1P_2$ 
where $\langle P_i,P_{i+1}\rangle \in (-1/2,1)$ for all $i\in \mathbb{Z}/3\mathbb{Z}$, and $\chi :=\langle P_i,P_{i+1}\times P_{i+2}\rangle >0$. From formula (\ref{formq2}) we see that  
$Q_2$ is directed towards the interior of $P_0P_1P_2$ precisely when the sign $\epsilon _2$ of 
$\langle Q_2,P_0\times P_1\rangle $ is positive. So $\epsilon _2=-1$ in Figure \ref{FIG1}. More generally, $Q_{i+2}$ is directed towards (away from) the interior of $P_0P_1P_2$ when the sign $\epsilon _{i+2}$ of $\langle Q_{i+2},P_i\times P_{i+1}\rangle$ is positive (negative).  

\spp
The centroid of the equilateral Euclidean triangle $P_0P_1Q_2$ is 
$$\frac{1}{3(1+\langle P_0,P_1\rangle))}((1+2\langle P_0,P_1\rangle )(P_0+P_1)+\epsilon _2\sqrt{1+2\langle P_0,P_1\rangle }P_0\times P_1).$$
So the equilateral spherical triangle $P_0P_1Q_2$ has centroid 
$$R_2:=\frac{(1+2\langle P_0,P_1\rangle )(P_0+P_1)+\epsilon _2\sqrt{1+2\langle P_0,P_1\rangle }P_0\times P_1}{\Vert  (1+2\langle P_0,P_1\rangle )(P_0+P_1)+\epsilon _2\sqrt{1+2\langle P_0,P_1\rangle }P_0\times P_1\Vert }=$$
$$\frac{\sqrt{1+2\langle P_0,P_1\rangle }(P_0+P_1)+\epsilon _2P_0\times P_1}{\sqrt{3}(1+\langle P_0,P_1\rangle )}.$$
Similarly, an equilateral spherical triangle $P_1,P_2Q_0$ has centroid 
$$R_0=\frac{\sqrt{1+2\langle P_1,P_2\rangle }(P_1+P_2)+\epsilon _0P_1\times P_2}{\sqrt{3}(1+\langle P_1,P_2\rangle )}.$$

\spp
So $3(1+\langle P_0,P_1\rangle )(1+\langle P_1,P_2\rangle )(1+\langle P_2,P_0\rangle )\langle R_2,R_0\rangle =(1+\langle P_2,P_0\rangle )S_1$ where  
$$S_1:= \langle \sqrt{1+2\langle P_0,P_1\rangle }(P_0+P_1)+\epsilon _2P_0\times P_1, \sqrt{1+2\langle P_1,P_2\rangle }(P_1+P_2)+\epsilon _0P_1\times P_2\rangle =$$
$$\sqrt{1+2\langle P_0,P_1\rangle }\sqrt{1+2\langle P_1,P_2\rangle }\langle P_0+P_1, P_1+P_2\rangle +$$
$$\epsilon _0 \sqrt{1+2\langle P_0,P_1\rangle }\langle P_0+P_1, P_1\times P_2\rangle +\epsilon _2 \sqrt{1+2\langle P_1,P_2\rangle }\langle P_0\times P_1, P_1+P_2\rangle +$$
$$ \epsilon _2\epsilon _0(1+\langle P_2,P_0\rangle )\langle P_0\times P_1, P_1\times P_2\rangle =$$
$$ \sqrt{1+2\langle P_0,P_1\rangle }\sqrt{1+2\langle P_1,P_2\rangle }(1+ \langle P_0,P_1\rangle +\langle P_2,P_0\rangle + \langle P_1,P_2\rangle ) 
+$$
$$ (\epsilon _0\sqrt{1+2\langle P_0,P_1\rangle }+\epsilon _2\sqrt{1+2\langle P_1,P_2\rangle } )\langle P_0, P_1\times P_2\rangle + \epsilon _2\epsilon _0(\langle P_0,P_1\rangle \langle P_1,P_2\rangle -\langle P_2,P_0\rangle)=$$
$$ \alpha \sqrt{1+2\langle P_0,P_1\rangle }\sqrt{1+2\langle P_1,P_2\rangle } 
+
 \chi  (\epsilon _0\sqrt{1+2\langle P_0,P_1\rangle }+\epsilon _2\sqrt{1+2\langle P_1,P_2\rangle } ) +
 \epsilon _2\epsilon _0(\langle P_0,P_1\rangle \langle P_1,P_2\rangle -\langle P_2,P_0\rangle ),$$
where 
$\alpha :=1+ \langle P_0,P_1\rangle +\langle P_2,P_0\rangle + \langle P_1,P_2\rangle $ is symmetric with respect to permutations of $P_0,P_1,P_2$, and the scalar triple product $\chi  :=\langle P_0, P_1\times P_2\rangle $ is symmetric with respect to cyclic permutations of $i=0,1,2\in \mathbb{Z}/3\mathbb{Z}$.

\spp 
Setting $d_i^2:=1+2\langle P_{i+1},P_{i+2}\rangle $ with $d_i\in (0,\sqrt{3})$, and  
$\gamma := 3(d_0^2+1)(d_1^2+1)(d_2^2+1)$, our calculation gives the useful formula
\begin{eqnarray}\label{useful}\frac{\gamma \langle R_{i+2},R_i\rangle}{d_{i+1}^2+1} &=& 4(\alpha d_{i+2}d_i  +
 \chi  (\epsilon _{i}d_{i+2}+\epsilon _{i+2}d_i ) )+\epsilon _{i+2}\epsilon _i((d_{i+2}^2-1)(d_i^2-1) -2(d_{i+1}^2-1)).\end{eqnarray}

\spp
By (\ref{eq20}),  
$$P_2= \frac{\langle P_2,P_0\rangle -\langle P_0,P_1\rangle \langle P_1,P_2\rangle }{1-\langle P_0,P_1\rangle ^2}P_0+\frac{\langle P_1,P_2\rangle-\langle P_0,P_1\rangle \langle P_2,P_0\rangle}{1-\langle P_0,P_1\rangle ^2}P_1+bP_0\times P_1$$
where by (\ref{eq1}), $b^2(1-\langle P_0,P_1\rangle ^2)^2=$
$$1-\langle P_0,P_1\rangle ^2-\langle P_1,P_2\rangle ^2-\langle P_2,P_0\rangle ^2+2\langle P_2,P_0\rangle \langle P_0,P_1\rangle \langle P_1,P_2\rangle .$$
Therefore, and because $\chi >0$, 
$$\chi  =b(1-\langle P_0,P_1\rangle ^2)=\sqrt{1-\langle P_0,P_1\rangle ^2-\langle P_1,P_2\rangle ^2-\langle P_2,P_0\rangle ^2+2\langle P_2,P_0\rangle \langle P_0,P_1\rangle \langle P_1,P_2\rangle }.$$
Summarising, we can say  
\begin{eqnarray}
\label{firstalphaeq} 2\alpha &=&d_0^2+d_1^2+d_2^2-1,\\
\label{betaeq} 4\chi  ^2&=&2(1-\alpha)(1+2\alpha)+d_0^2d_1^2+d_1^2d_2^2+d_2^2d_0^2+d_0^2d_1^2d_2^2, 
\end{eqnarray}

 \section{Proof of Theorem \ref{THEOREM1}}\label{SECTION:3} 
Note that $P_0P_1P_2$ is equilateral precisely when $d_0=d_1=d_2$. This case is excluded in the present section, namely we suppose that $d:=(d_0,d_1,d_2)\in D:=(0,\sqrt{3})^3\setminus  \Delta $, 
where $\Delta:=\{ (\delta ,\delta ,\delta):\delta \in (0,\sqrt{3})\}$. 
Set\footnote{Because of these choices, if $\langle P_2,P_0\rangle =\langle P_1,P_0\rangle$ then $\langle R_2,R_0\rangle =\langle R_1,R_0\rangle$. Similarly, if $\langle P_0,P_1\rangle =\langle P_2,P_1\rangle$ then $\langle R_0,R_1\rangle =\langle R_2,R_1\rangle$, and if $\langle P_1,P_2\rangle =\langle P_0,P_2\rangle$ then $\langle R_1,R_2\rangle =\langle R_0,R_2\rangle$. At this stage however, our attention is not limited to the case where $P_0P_1P_2$ is isosceles.} $\epsilon _0=\epsilon _1=\epsilon _2=\epsilon :=\pm1$. So the equilateral triangles constructed on the edges of $P_0P_1P_2$ are directed either all outwards ($\epsilon =-1$) or all inwards ($\epsilon =+1$).  By (\ref{useful}),  
\begin{eqnarray}\label{r2r0}\frac{\gamma \langle R_{i+2},R_i\rangle }{d_{i+1}^2+1}&=&  4(\alpha d_{i+2}d_i +
\epsilon \chi  (d_{i+2}+d_i)) +(d_{i+2}^2-1)(d_i^2-1) -2(d_{i+1}^2-1).\end{eqnarray}
A short calculation then shows that $\langle R_2,R_0\rangle =\langle R_0,R_1\rangle$ precisely when 
$$ 4\left( \alpha d_0(1-d_1d_2) +\epsilon \chi   (1-d_0d_1-d_1d_2-d_2d_0)\right) (d_2-d_1)
 =  
-2 (d_2^2-d_1^2)(d_0^2+d_1^2+d_2^2-1). $$
The right hand side is $-4\alpha  (d_2^2-d_1^2)$, and so our condition is that either $d_1=d_2$ or 
\begin{eqnarray}\label{r0r1eqg}  \alpha (d_0+d_1+d_2-d_0d_1d_2) +\epsilon \chi   (1-d_0d_1-d_1d_2-d_2d_0)
 &=& 0. \end{eqnarray}
%
%
 %
 This equation is invariant with respect to cyclic permutations of the subscripts and, because $d\notin \Delta$, equation (\ref{r0r1eqg}) is necessary and sufficient for $R_0R_1R_2$ to be equilateral. 
We observe that:
\begin{enumerate}
\item For $d\in D$ we have $d_0+d_1+d_2-d_0d_1d_2>0$ because otherwise, by the AM-GM inequality, $$d_0d_1d_2\geq d_0+d_1+d_2\geq 3(d_0d_1d_2)^{1/3} \Longrightarrow d_0d_1d_2\geq 3\sqrt{3},$$
contradicting $d_0,d_1,d_2\in (0,\sqrt{3})$. 
\item For any given $d\in D$, equation (\ref{r0r1eqg}) cannot hold for both $\epsilon =1$ and $\epsilon =-1$. Otherwise, 
$$ \alpha (d_0+d_1+d_2-d_0d_1d_2)=0= \chi   (1-d_0d_1-d_1d_2-d_2d_0)\Longrightarrow$$
$$d_0^2+d_1^2+d_2^2=1=d_0d_1+d_1d_2+d_2d_0\Longrightarrow d_0+d_1+d_2=\sqrt{3}\Longrightarrow $$
$$\langle d,(1,1,1)\rangle ^2=3=\langle d,d\rangle \langle (1,1,1),(1,1,1)\rangle \Longrightarrow d=(1,1,1)/\sqrt{3}$$
by Cauchy-Schwarz, contradicting $d\notin \Delta$. 
\item So (\ref{r0r1eqg}) holds neither for any $d\in D\cap S^2_+$, nor for any $d\in D$ satisfying
$$d_0d_1+d_1d_2+d_2d_0=1.$$
\item $\epsilon$ has the sign of 
$$(1-d_0^2-d_1^2-d_2^2)(1-d_0d_1-d_1d_2-d_2d_0).$$
\end{enumerate}
Hence, a necessary condition for $P_0P_1P_2$ to be either outward-Napoleonic or inward-Napoleonic is
\begin{equation}\label{algeq}
\alpha ^2(d_0+d_1+d_2-d_0d_1d_2)^2-\chi  ^2 (1-d_0d_1-d_1d_2-d_2d_0)^2=0
\end{equation} 
with $\alpha ,\chi ^2$ given by (\ref{firstalphaeq}), (\ref{betaeq}) as symmetric polynomials in $d_0,d_1,d_2$. 

\spp
Substituting for $\alpha,\chi ^2$, a calculation shows that the left-hand side of (\ref{algeq}) is $\gamma /12$ times 
\begin{equation}\label{factalgeq}
(d_0^2+d_1^2+d_2^2-d_0d_1-d_1d_2-d_2d_0)(d_0^2+d_1^2+d_2^2+d_0d_1+d_1d_2+d_2d_0-2),
\end{equation}
where the first factor is nonzero because $d\notin \Delta$. So $P_0P_1P_2$ is either outward-Napoleonic or inward-Napoleonic when
$$(d_0+\frac{d_1}{2}+\frac{d_2}{2})^2+\frac{3}{4}(d_1+\frac{d_2}{3})^2 +\frac{2}{3}d_2^2=2.$$
When this happens we find that, using 
$d_0^2+d_1^2+d_2^2-1=1-d_0d_1-d_1d_2-d_2d_0$,
$$(1-d_0^2-d_1^2-d_2^2)(1-d_0d_1-d_1d_2-d_2d_0)=-(1-d_0^2-d_1^2-d_2^2)^2<0$$
by observations (2), (3). Then $\epsilon =-1$ by observation (4), namely $P_0P_1P_2$ is outward-Napoleonic (it cannot be inward-Napoleonic).   

\spp
More can be said: using $1-d_0d_1-d_1d_2-d_2d_0=d_0^2+d_1^2+d_2^2-1=2\alpha$ and $\epsilon =-1$ in (\ref{r0r1eqg}), we find that
$$2\chi =d_0+d_1+d_2-d_0d_1d_2.$$  
\spp
Then, substituting for $\chi ,\alpha$, we find  
$$(d_0^2+1)(d_2^2+1)+(4(\alpha d_2d_0-\chi (d_2+d_0))+(d_2^2-1)(d_0^2-1)-2(d_1^2-1)=$$
$$2(d_0d_2-1)(d_0^2+d_1^2+d_2^2+d_0d_1+d_1d_2+d_2d_0-2)=0.$$
Comparing the left-hand side with equation (\ref{r2r0}) for $i=0$,   
$$\langle R_2,R_0\rangle =\langle R_0,R_1\rangle =\langle R_1,R_2\rangle =-1/3.$$
So the Napoleonisation $R_0R_1R_2$ belongs to the congruency class of spherical equilateral triangles of perimeter
$3(\pi -\arccos (1/3))$. This proves Theorem \ref{THEOREM1}.~$\square$

\section{Context}\label{contextsec}
\spp
The study of triangular patterns arises in several situations and we present here
some of them, both in the Euclidean and spherical case, to provide context to the
Napoleonic construction.

\spp
To start with, the theory of crystalline growth provides a number of examples of triangular structures,
which occur often in nature, especially for mineral crystals. These structures can arise as sections of trigonal crystal shapes in ruby, quartz, brucite, calcite, agate, jasper, etc., and in groups of sodium nitrate and nickel arsenide, see~\cite[page~11, page~46, Figure~A.4]{SSC}.
They also occur as sections of tetragonal systems of silicon, see~\cite{PAUL}, and in hexagonal patterns of beryl, cancrinite, apatite, sugilite, etc., see~\cite{ATLAS}, and of oxygen
atoms in hexagonal	lattices of ice, see~\cite[Figure~2.3]{PETRE} and~\cite[pages 4-5]{DEER}. 

\spp
With regards to organic molecules, the case of the compound theophylline is also noteworthy, since, differently from most organic compounds, it may form crystals with a triangular form, see~\cite{UNUS} (and in particular Figure~1 there for optical microscope images showing the growth of triangular theophylline crystals on a glass slide). In terms of applications, theophylline is a bronchodilator and is used as a medication to treat asthma and obstructive pulmonary diseases. 
Triangular crystals also occur in the crystalline growth of multilayers silver nanocrystals, see~\cite{Courty2007}.

\spp
Another interesting situation arises in the theory of superfluids and superconductors, in which
vortices have been found to form triangular lattice arrangements, as a configuration that minimises the total energy, see~\cite{zbMATH03188431, BEN, CORR}
(remarkably, in this setting, the original idea of~\cite{Abrikosov}
that the energy is minimised by a square lattice turns out to be inaccurate).

\begin{figure}[h]
\centering
\includegraphics[height=3.7cm]{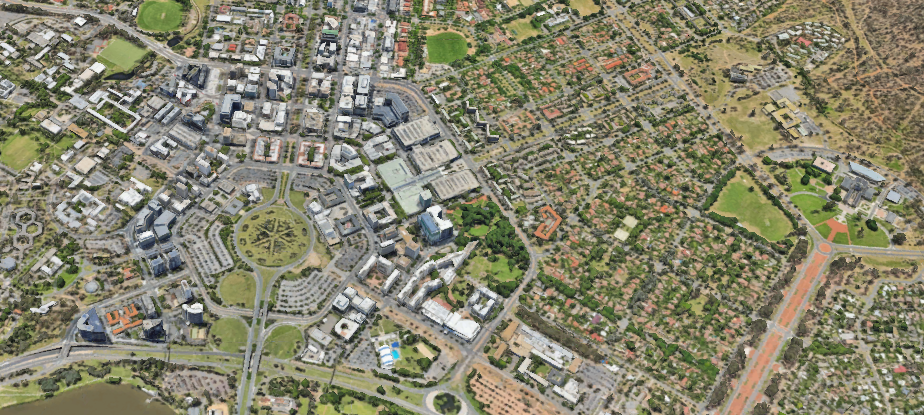}
\includegraphics[height=3.7cm]{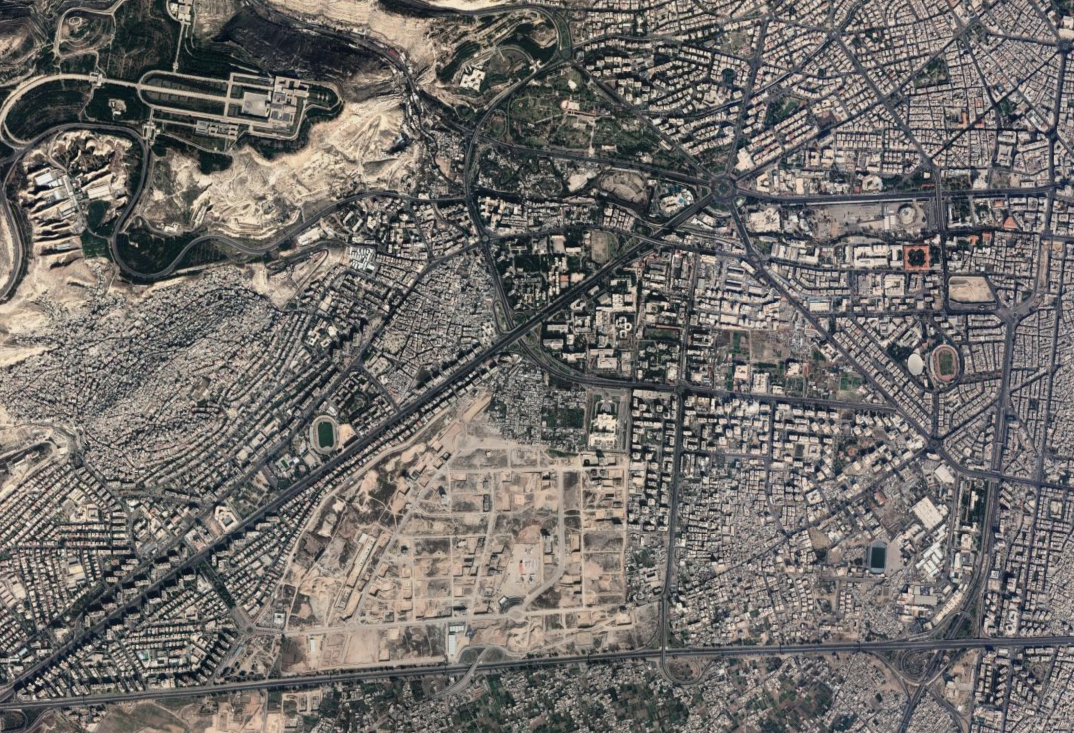} 
\includegraphics[height=3.7cm]{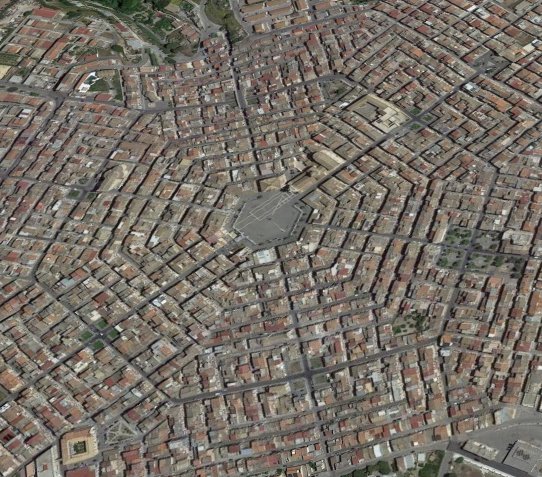} \vskip2pt
\includegraphics[height=3.7cm]{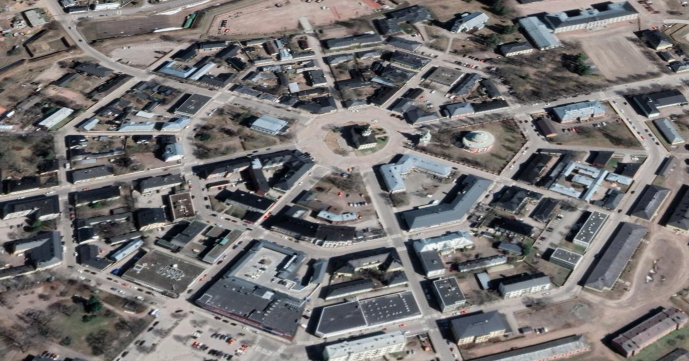}
\includegraphics[height=3.7cm]{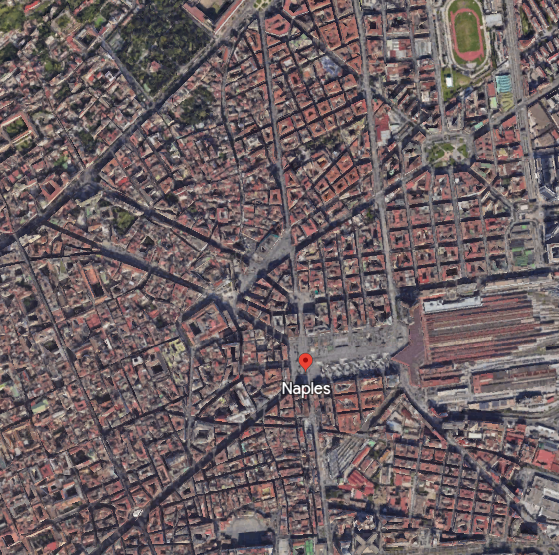} 
\includegraphics[height=3.7cm]{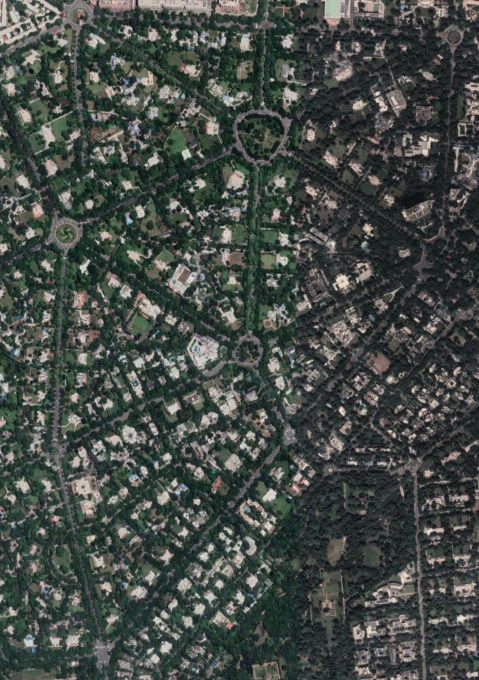}
\includegraphics[height=3.7cm]{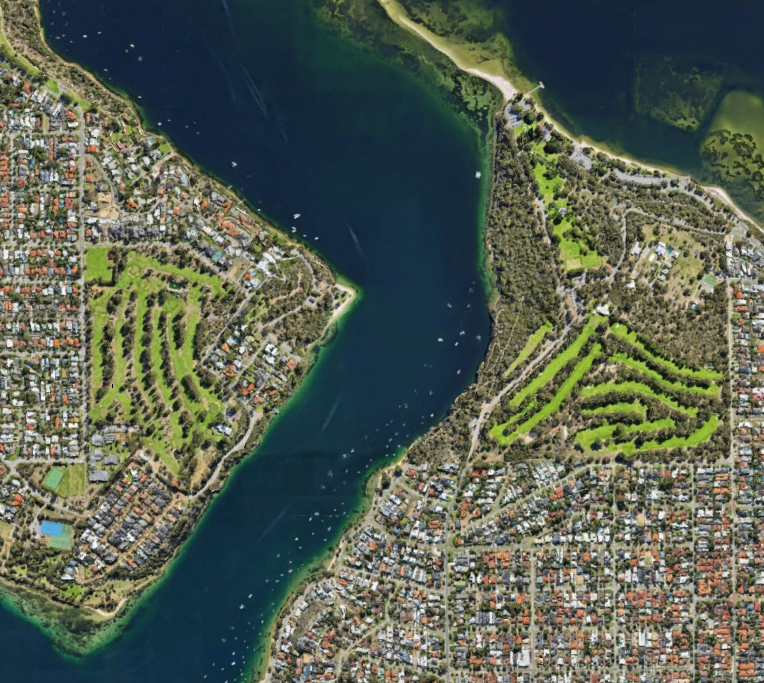}\vskip2pt 
\includegraphics[height=3.7cm]{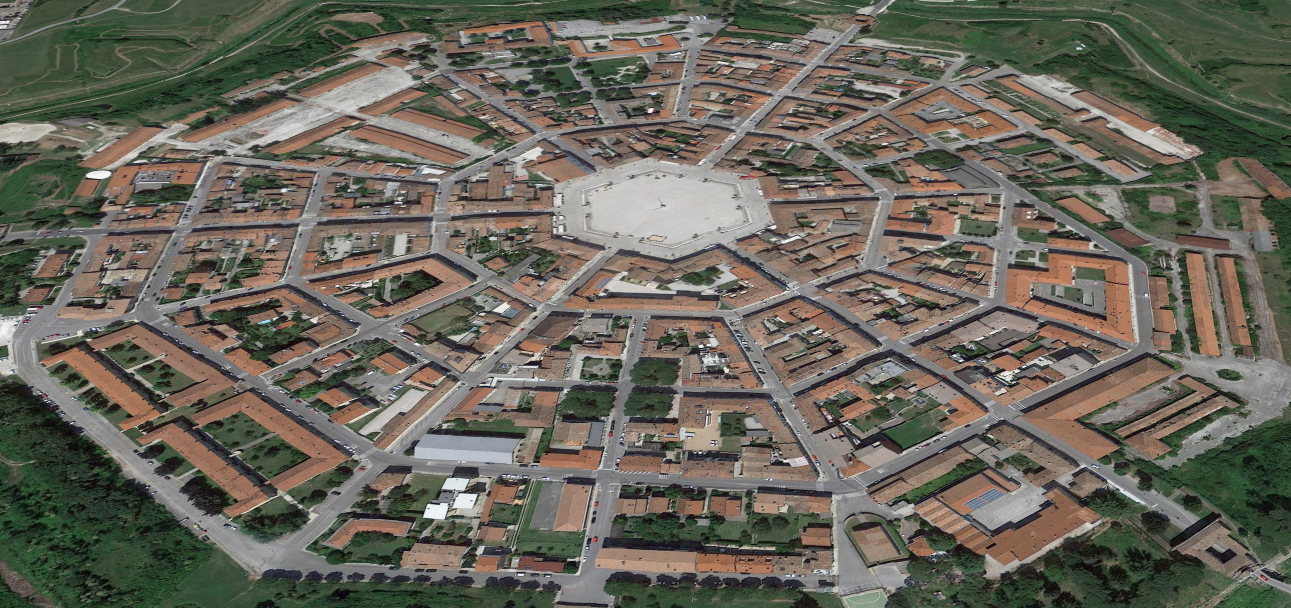} 
\includegraphics[height=3.7cm]{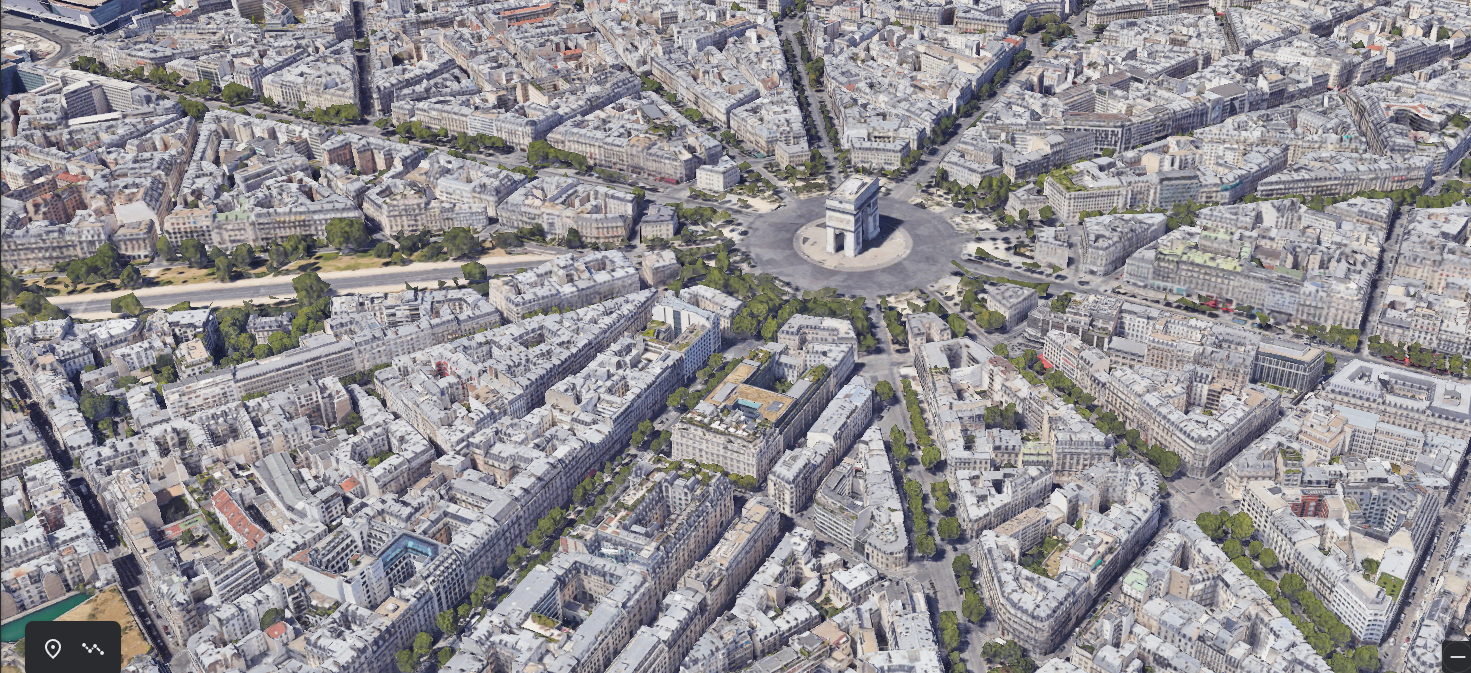}\vskip2pt
\includegraphics[height=3.7cm]{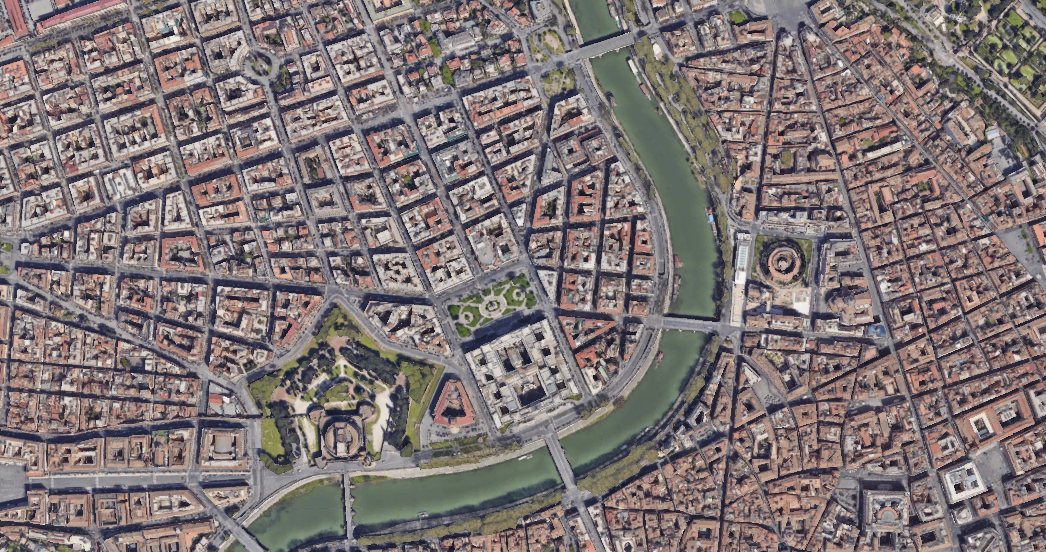}
\includegraphics[height=3.7cm]{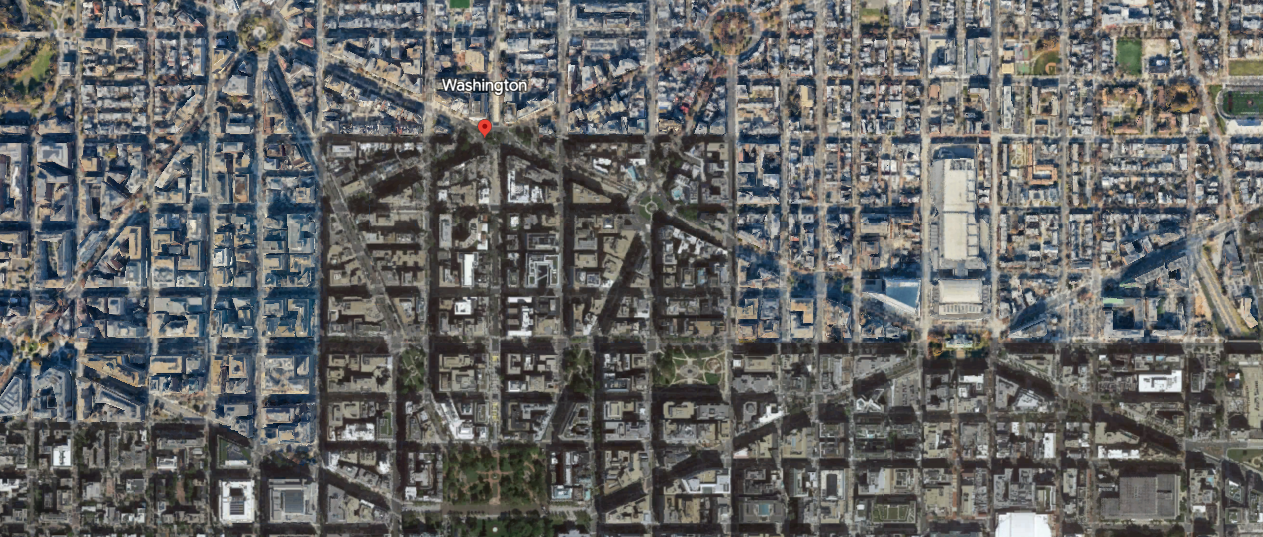}
\caption{\em Aerial views of
Canberra, Damascus, Gramichelle, Hamina, Naples, New Dheli, Perth, Palmanova, Paris, Rome, Washington
(screenshots from Google Earth).}
\label{3stPERINDETRAPf-A}\end{figure}

\begin{figure}[h]
\centering
\includegraphics[height=7.1cm]{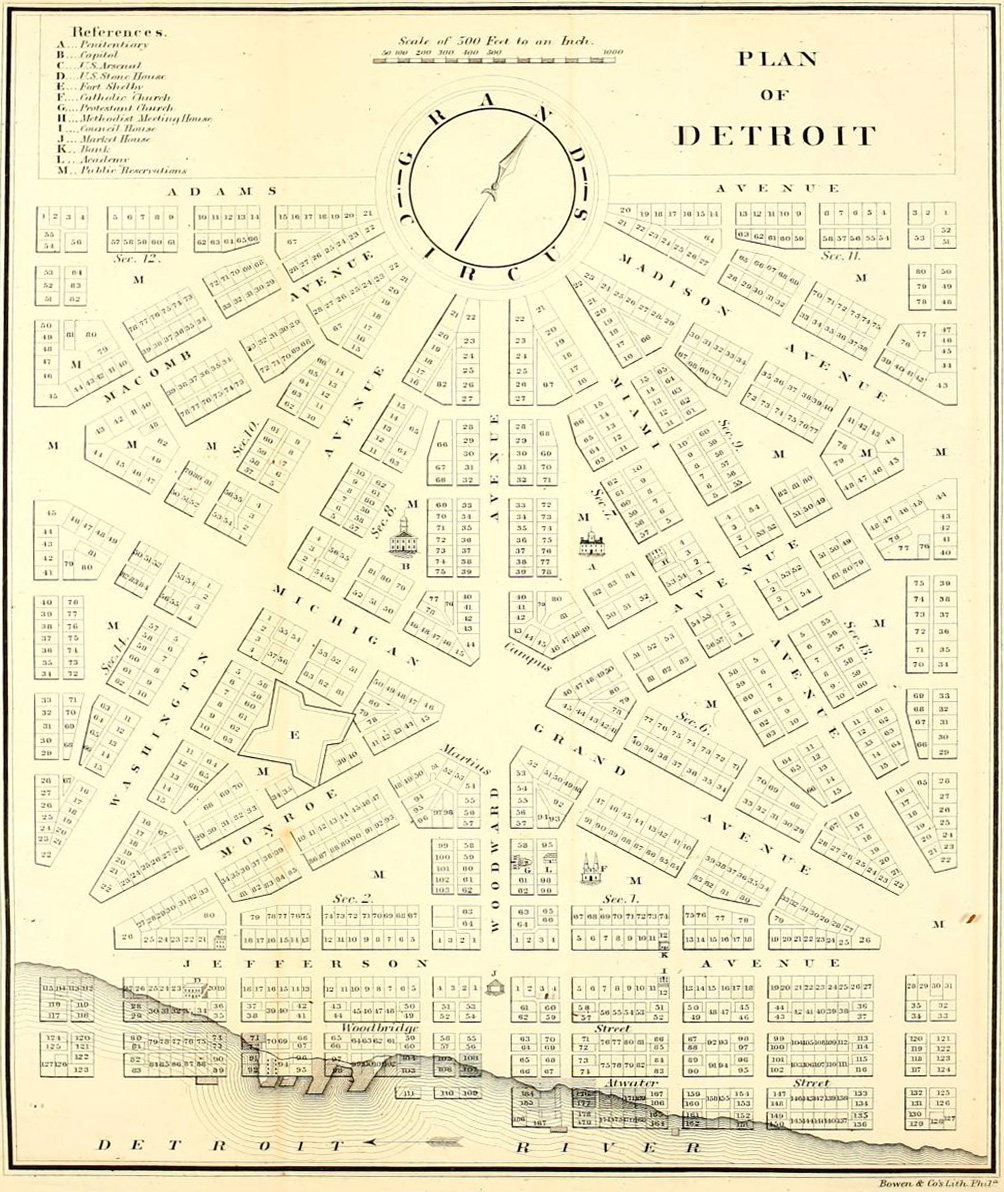}
\caption{\em Detroit city layout plan following the 1805 fire that destroyed most of the city
(public domain image from Wikipedia).}
\label{3stPERINDETRAPf}\end{figure}

\spp
Further context comes in urban design and town planning,
since several cities exhibit triangle-shaped blocks,
see Figures~\ref{3stPERINDETRAPf-A} and~\ref{3stPERINDETRAPf}. We refer to~\cite{HES} for more information about hexagonal and triangular
plannings in urban design.
A nice scenario in which Napoleonic constructions arise in the context of city planning is therefore the case of an initial block given by a (possibly irregular) triangle, which gets developed by three adjacent blocks formed by equilateral triangles. If, for instance, we assume that each of these new blocks is endowed with a communication antenna located at the centre of the block, the fundamental Napoleonic question is whether or not these antennas form an equilateral\footnote{Indeed, one can suppose that a rewarding strategy in signal transmission consists in ``covering
the largest possible surface with the least possible mutual distances among the antennas'':
the formulation of this isoperimetric problem in the class of triangular patterns of antennas is optimised by
equilateral configurations (in turn, the isoperimetric problem for triangles can
be explicitly solved, for example, by writing the area of a general triangle using Heron's formula
and using the AM-GM inequality to see that the optimiser is equilateral;
another proof can be performed by purely geometrical arguments by considering an ellipse with foci 
centred at two vertices and passing through the third one).} triangle.

\begin{figure}[h]
\centering
\includegraphics[height=6.1cm]{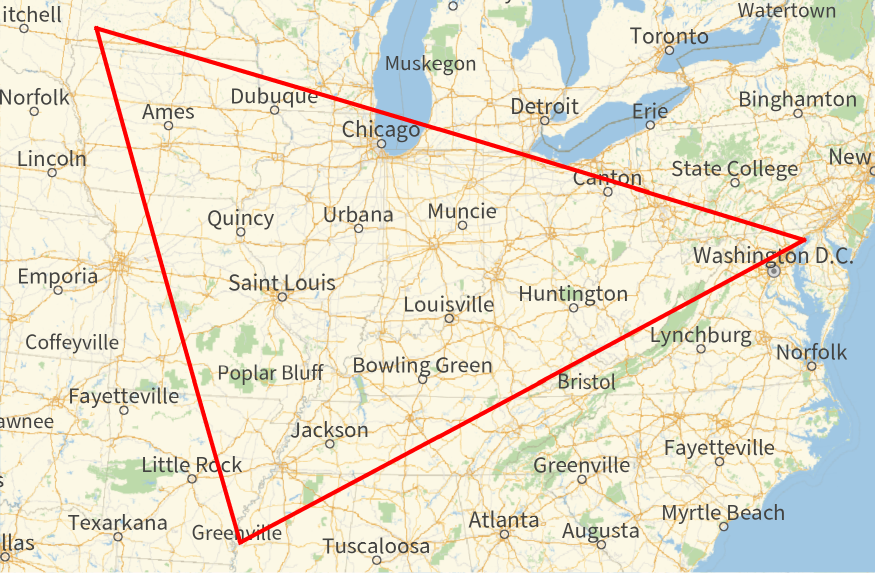}$\qquad$
\includegraphics[height=6.1cm]{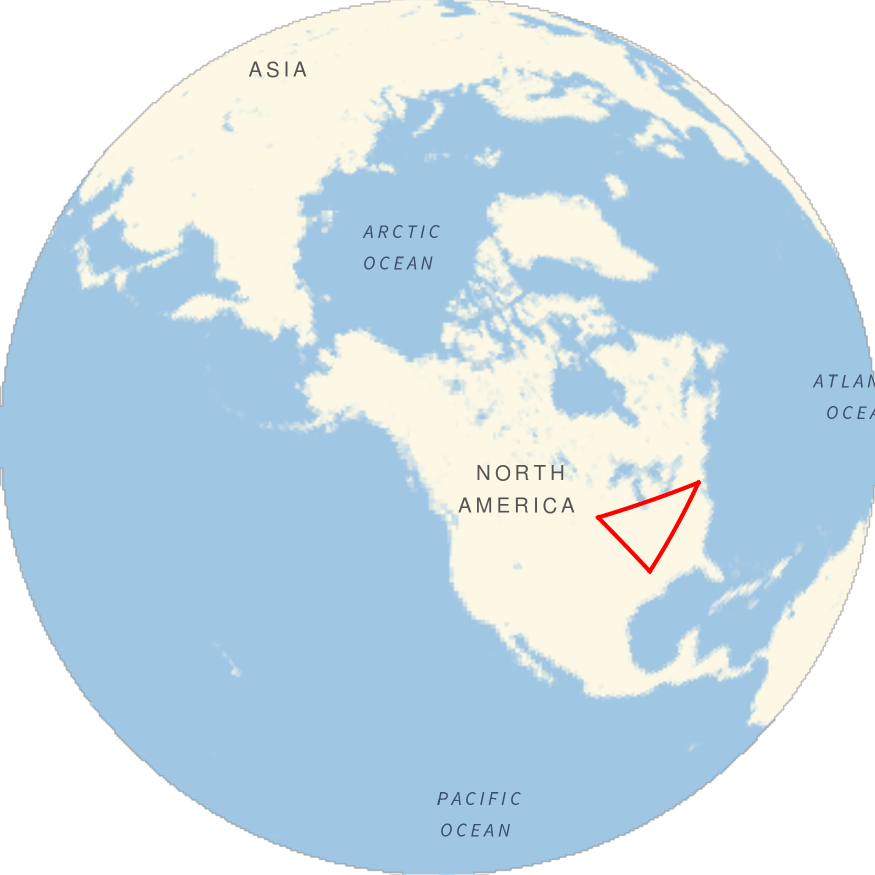}
\caption{\em Location of Aberdeen Proving Ground, Maryland, Chandler Air Force Station, Minnesota, and Greenville Air Force Base, Mississippi (obtained via Mathematica).}
\label{3stPERINAPf}\end{figure}

\spp
Besides the interest in mathematical theories, it is worth recalling that
triangular patterns of antennas are broadly used in modern technologies. For instance,
a triangulation network is implemented in
North America for satellite communication, and
satellite geodesy exploited an equilateral triangle on the Earth with operational basis established in 1963
at Aberdeen Proving Ground
39$^\circ$28'24.4''N 76$^\circ$08'27.0''W, Maryland, Chandler Air Force Station
43$^\circ$53'52''N 95$^\circ$56'45''W, Minnesota, and Greenville Air Force Base 33$^\circ$28'58''N 90$^\circ$59'08''W, Mississippi: these stations form
a triangle with about 1500 km on a side, see Appendix~E of the North American
Datum report {\url{https://www.ngs.noaa.gov/PUBS_LIB/NAD_National_Academy_of_Science.pdf}} and Figure~\ref{3stPERINAPf}
(actually, the geodesic distance between 
39$^\circ$28'24.4''N 76$^\circ$08'27.0''W and
43$^\circ$53'52''N 95$^\circ$56'45''W is about 1700 km,
the one between 43$^\circ$53'52''N 95$^\circ$56'45''W and
33$^\circ$28'58''N 90$^\circ$59'08''W is about 1250 km,
and the one between
33$^\circ$28'58''N 90$^\circ$59'08''W
and
39$^\circ$28'24.4''N 76$^\circ$08'27.0''W
is about 1500 km).
Also, in the Laser Interferometer Space Antenna
(LISA), three spacecraft are arranged in an equilateral triangle, see~\cite{LISA}.
See also~\cite{MYSTE} for other different applications of equilateral triangles.

\spp
It is interesting to point out that
the Napoleonic construction is also related to the Fermat-Steiner problem, see~\cite[\S7.1]{MR2022170}. This problem was posed in 1629 by Fermat and addressed by Torricelli sometime before 1640. The question is to find a point whose total sum of distances from the vertices of a given triangle is a minimum. That is, given a triangle with vertices~$A$, $B$ and~$C$, find a point~$ X$ for which~$\|X-A\|+ \|X-B\|+ \|X-C\|$ is a minimum.
Remarkably, the point~$X$ can be obtained as the outcome of a Napoleonic construction. Indeed, one constructs the three outward pointing equilateral triangles on the sides of the triangle with vertices~$A$, $B$ and~$C$ and connects each vertex of the original triangle to the new vertex
of the opposite equilateral triangle, the common intersection of 
these segments thus
identifying the desired point~$X$
(this, if all the angles of the original triangle are less than~$\frac{2\pi}3$, 
otherwise~$X$ coincides with the vertex corresponding to the largest angle).
The point~$X$ found in this way is sometimes called ``Fermat point'' (or ``Torricelli point'', or ``Steiner point'').

\spp
The Fermat-Steiner problem itself is a source of questions of practical importance, such as the so-called ``factory problem'': suppose that there are three factories, labeled $A$, $B$, and~$C$, which need to be supplied by shipments of a certain material originating at a warehouse~$X$, and assume additionally that the cost of the shipment is proportional to the distance between the warehouse and the factory receiving the shipment. Once again, the total cost would be thus proportional to~$\|X-A\|+ \|X-B\|+ \|X-C\|$
and therefore the optimal location for~$X$ is that of the Fermat point above.

\spp
Equilateral triangular patterns naturally emerge
from chaotic processes (see~\cite[Figure~6.2]{CHF}). Besides,
other interesting occurrences of equilateral triangles arise for the Kakeya needle problem: namely,
equilateral triangles are the convex figures of least
area in which a line segment of given length 1 can be moved around and return
to its original position with reversed orientation, see~\cite[page~221]{GARD} and the references therein
(and Kakeya problems are nowadays mainstream in the contemporary mathematics research, see e.g.~\cite{MR2206765, MR3946717}).

\spp
Equilateral triangles also play a pivotal role in numerics and finite element methods, since, for numerical accuracy and efficiency, the shape of the finite elements is often chosen to resemble as much as possible the ``reference elements'', such as equilateral triangles or tetrahedra, see e.g.~\cite[pages~106 and~137]{CAD} and~\cite[page~242]{KURO}.

\spp
Extensions of the Napoleonic construction to polygons (rather than triangles) are also possible.
The most famous result in this framework is the so-called PDN-Theorem, named after
Karel Petr, Jesse Douglas, and Bernhard Neumann ~\cite{PETR-zbMATH02641104, DOUG-MR2178, NEUM-MR6839} (the special case of quadrilaterals being related to Van Aubel's Theorem, see~\cite{VANB}, and
to Th\'{e}bault's Second Problem, see~\cite{MR3961650}).

\spp
In all these settings, the case of the flat plane is somewhat special and it is intriguing
to understand the role played by a possibly curved environment: the case of constant curvature being the natural
one to start with, we study here spherical Napoleonic triangles (we also aim at discussing hyperbolic spaces in a forthcoming project).


\end{document}